\documentclass{amsart}
\usepackage{amssymb, amsmath, amsfonts, amscd}

\input xypic

\theoremstyle{plain}
\newtheorem{theorem}{Theorem}[section]
\newtheorem{corollary}[theorem]{Corollary}
\newtheorem{lemma}[theorem]{Lemma}
\newtheorem{proposition}[theorem]{Proposition}

\newtheorem{example}[theorem]{Example}

\theoremstyle{definition}
\newtheorem{definition}[theorem]{Definition}

\theoremstyle{remark}

\numberwithin{equation}{theorem}

\newcommand{\m}{\mathfrak{m}}
\newcommand{\F}{\mathcal{F}}

\newcommand{\I}{\mathcal{I}}

\renewcommand{\L}{\mathcal{L}}
\newcommand{\E}{\mathcal{E}}

\newcommand{\Q}{\mathcal{Q}} 
\renewcommand{\O}{\mathcal{O} }
\renewcommand{\a}{\alpha}

\renewcommand{\Pr}{\mathcal{P} }

\newcommand{\SL}{\operatorname{SL}} 
\newcommand{\GL}{\operatorname{GL}} 
\newcommand{\Pic}{\operatorname{Pic} }

\newcommand{\End}{\operatorname{End} }
\newcommand{\Spec}{\operatorname{Spec} }

\renewcommand{\H}{\operatorname{H} }
\newcommand{\Proj}{\operatorname{Proj} }
\newcommand{\D}{\operatorname{D} }

\newcommand{\R}{\operatorname{R} }
\newcommand{\U}{\operatorname{U}}

\renewcommand{\lg}{\mathfrak{g}}
\newcommand{\lp}{\mathfrak{p}} 
\renewcommand{\a}{\mathfrak{a}}

\renewcommand{\ln}{\mathfrak{n}}
\renewcommand{\l}{\mathfrak{l}}
\newcommand{\lh}{\mathfrak{h}}

\newcommand{\nplus}{\mathfrak{n}_{+}}
\newcommand{\nminus}{\mathfrak{n}_{-}}

\newcommand{\onep}{\mathbf{1}_{\lp}}
\newcommand{\oneg}{\mathbf{1}_{\lg}}

\newcommand{\gr}{\mathbb{G} }
\newcommand{\sym}{\operatorname{Sym} }

\newcommand{\Vl}{V_\lambda}

\newcommand{\ogpmod}{\underline{mod}$^G(\O_{G/P}$) }
\renewcommand{\pmod}{\underline{mod}($P$) }
\newcommand{\oghmod}{\underline{mod}$^G(\O_{G/H}$) }
\newcommand{\hmod}{\underline{mod}($H$) }

\newcommand{\p}{\mathbb{P}}

\begin{document}

\title{On jet bundles and generalized Verma modules}

\author{Helge Maakestad}
\address{Institut Fourier, Grenoble}

\email{maakestd@fourier.ujf-grenoble.fr}


\date{December 2008}

\begin{abstract} Let $K$ be a field of characteristic zero and let
  $W\subseteq V$ be $K$-vector spaces of dimension $m$ and $m+n$. Let
  $P\subseteq \SL(V)=G$ be the subgroup fixing $W$. It follows
  $X=G/P$ equals the grassmannian of $m$-planes in $V$. Let $\O_X(d)$
  be the $d$'th tensor power of the line bundle $\O(1)$ coming from
  the Pl\"{u}cker embedding of $X$. There is an
  equivalence of categories between the category of finite dimensional
  $P$-modules and the category of $G$-linearized locally free finite rank
  $\O_{X}$-modules. The $l$'th jet bundle $\Pr^l_X(\O_X(d))$ 
  is a $G$-linearized locally free sheaf and the aim of this paper is to
  describe its corresponding $P$-module using higher direct images of
  $G$-linearized sheaves, filtrations of generalized Verma modules,
  canonical filtrations of irreducible $\SL(V)$-modules and annihilator ideals of
  highest weight vectors and to apply this to the study of
  discriminants of linear systems on grassmannians.
The main theorem of the paper is that the discriminant $D^l(\O_X(d))$
is irreducible for all $1\leq l < d$.
\end{abstract}

\maketitle

\tableofcontents

\section{Introduction} 

Let $K$ be a field of characteristic zero and let
  $W\subseteq V$ be $K$-vector spaces of dimension $m$ and $m+n$. Let
  $P\subseteq \SL(V)=G$ be the subgroup fixing $W$. It follows
  $X=G/P$ equals the grassmannian of $m$-planes in $V$. There is an
  equivalence of categories between the category of finite dimensional
  $P$-modules and the category of $G$-linearized locally free finite rank
  $\O_{X}$-modules. The $l$'th jet bundle $\Pr^l_X(\L)$ where
  $\L\in \Pic^{G}(X)$ is a $G$-linearized locally free $\O_X$-module, and the aim of this paper is to
  describe its corresponding $P$-module using Taylor morphisms, higher direct images of
$G$-linearized sheaves, filtrations of generalized Verma modules,
canonical filtrations of irreducible $\SL(V)$-modules and annihilator ideals of
highest weight vectors. 

Let $i:X\rightarrow \p(\wedge^m V^*)$ be the
Pl\"{u}cker embedding and let $\O_X(d)=i^*\O(1)^{\otimes d}$ where $\O(1)$
is the tautological line bundle on $\p(\wedge^m V^*)$.
In Theorem \ref{mainp} we prove there is an isomorphism
\[ \Pr^l_X(\O_X(d))(x)^*\cong \U_l(\lg)v \]
of $P$-modules. Here $\Pr^l_X(\O_X(d))(x)^*$ is the dual of the fiber
of the jet bundle at the distinguished point $x\in G/P$  and
$\U_l(\lg)v\subseteq \H^0(X,\O_X(d))^*$ is the $l$'th piece of the
canonical filtration of $\H^0(X,\O_X(d))^*$.
We then apply the results on jet bundles and Taylor maps to the study of discriminants of
line bundles on grassmannians. In Theorem \ref{irreducible} we prove
the $l$'th discriminant $D^l(\O_X(d))$ is irreducible when $1\leq l
\leq d$.

The paper is organized as follows:
In section two of the paper we study general properties of jet bundles
of $G$-linearized locally free sheaves on homogeneous spaces
$G/H$. Here $G$ is a linear algebraic group of finite type over a
field $K$ and $H\subseteq G$ is a closed subgroup.

In section three of the paper we study the \emph{Taylor morphism}
$T^l$ for invertible sheaves on the grassmannian.
Let  $\O_X(d)$ be the $d$'th tensor power of the tautological bundle
coming from the Pl\"{u}cker embedding. 
We prove in Corollary \ref{surjective} the Taylor morphism 
\[ T^l:\H^0(X,\O_X(d))\otimes \O_X\rightarrow \Pr^l_X(\O_X(d)) \]
is a surjective map of locally free sheaves when $1\leq l \leq d$.

In section four we study the \emph{canonical filtration} $\U_l(\lg)v$ of
of the irreducible $G$-module $\H^0(X,\O_X(d))^*$. Using
the universal enveloping algebra
$\U(\lg)$ and the annihilator ideal $ann(v)\subseteq \U(\lg)$ where $v\in \H^0(X,\O_X(d))^*$
is the highest weight vector we give in Corollary \ref{dim} a basis for $\U_l(\lg)v$ as
$K$-vector space. We also compute the dimension of $\U_l(\lg)v$.

In section five we study the dual of the fiber of the jet bundle
$\Pr^l_X(\O_X(d))(x)^*$ at the distinguished point $x\in X=G/P$ as $P$-module. 
Using the results obtained in the previous sections, we classify in
Theorem \ref{mainp} the fiber as $P$-module, and describe it in terms
of $\U_l(\lg)v$ - a subquotient of the generalized Verma module $\U(\lg)\otimes L_v$
associated to the line $L_v$ spanned by the highest weight vector $v\in \H^0(X,\O_X(d))^*$.

In section six we apply the results obtained in the previous sections
and study discriminants of line bundle on the grassmannian $X$. 
We prove the $l$'th discriminant $D^l(\O_X(d))$ is irreducible
when $1\leq l \leq d$.

The motivation for the study of the $P$-module $\Pr^l_X(\O_X(d))(x)^*$
is a relationship between the jet bundle and resolutions of ideal sheaves
of discriminants of linear systems  (see \cite{maa10} Example 5.12). Using the jet
bundle $\Pr^l_X(\O_X(d))$ one constructs a double complex
\[ \O_X(-j)\otimes \H^i(X,\wedge^j\Pr^l_X(\O_X(d))^*) \]
of sheaves on $\p(W^*)$ where $W=\H^0(X,\O_X(d))$.
The $l$'th discriminant
$D^l(\O_X(d))$ of the line bundle $\O_X(d)$ is a closed subscheme
\[ D^l(\O_X(d))\subseteq \p(W^*) \]
and the double complex $\O_X(-j)\otimes
\H^i(X,\wedge^j\Pr^l_X(\O_X(d))^*) $ may in some cases be used to construct a resolution
of the ideal sheaf of $D^l(\O_X(d))$. Knowledge on the $P$-module
structure of $\Pr^l_X(\O_X(d))(x)^*$ will give information on the
problem of constructing such a resolution.

\section{Jet bundles on quotients}

In this section we study general properties of jet bundles on
quotients $G/H$. Here $G$ is a linear algebraic group of finite type
over a field and $H\subseteq G$ a closed subgroup.

\textbf{Notation:} Let $\pi:X\rightarrow S$ be a smooth
and separated morphism. Let $Y=X\times_S X$ and let 
\[ \Delta:X\rightarrow Y \]
be the diagonal closed embedding. Let $p,q:Y\rightarrow X$ be the
canonical projection maps and let $\I_\Delta \subseteq \O_Y$ be the
ideal of the diagonal. Let $\O_{\Delta^l}=\O_Y/\I_\Delta^{l+1}$ the
structure sheaf of the \emph{$l$'th infinitesimal neigborhood of the diagonal}.

There is on $Y$ an exact sequence of $\O_Y$-modules
\begin{align}
&\label{short} 0\rightarrow \I_{\Delta}^{l+1}\rightarrow \O_Y
\rightarrow \O_{\Delta^{l+1}} \rightarrow 0.
\end{align}
Applying the functors $\R^ip_*(-\otimes q^*\E)$ and the formalism of
derived functors to \ref{short} we get a long exact sequence of quasi coherent $\O_X$-modules
\begin{align}
&\label{longexact} 0\rightarrow p_*(\I_{\Delta}^{l+1}\otimes q^*\E)\rightarrow
p_*q^*\E \rightarrow p_*(\O_{\Delta^{l+1}}\otimes q^*\E) \rightarrow
\end{align}
\[\R^1p_*(\I_{\Delta}^{l+1}\otimes q^*\E)\rightarrow \R^1p_*q^*\E
\rightarrow \R^1p_*(\O_{\Delta^{l+1}}\otimes q^*\E) \rightarrow \cdots
\]
By flat basechange it follows
\[ \R^ip_*q^*\E\cong \pi^*\R^i\pi_*\E. \]
for all $i\geq 0$.
It follows 
\[ p_*q^*\E\cong \pi^*\pi_*\E\otimes \O_X.\]

\begin{definition} \label{jets} The quasi coherent $\O_X$-module
\[ \Pr^l_{X/S}(\E)=p_*(\O_{\Delta^{l+1}}\otimes q^*\E) \]
is the \emph{$l$'th order jet bundle} of $\E$. The morphism
\[ T^l:\pi^*\pi_*\E\otimes \O_X\rightarrow \Pr^l_{X/S}(\E) \]
is the \emph{$l$'th Taylor morphism} for $\E$.
\end{definition}

Note: the $\O_X$-module $\Pr^l_{X/S}(\E)$ has a left and right
$\O_X$-structure. In this paper we will only consider the left $\O_X$-structure.

\begin{example} Taylor maps on smooth projective schemes.\end{example}

Let $S=\Spec(K)$ and  $X\subseteq \p^n_S$ be a smooth projective scheme of finite type
over the field $K$. Let $\E$ be a locally free $\O_X$-module. 
The Taylor morphism for this situation looks as follows:
\[ T^l :\H^0(X,\E)\otimes \O_X \rightarrow \Pr^l_{X/S}(\E) .\]
Given a $K$-rational point $x\in X$ and a global section $s\in
\H^0(X,\E)$ it follows the Taylor map $T^l$ formally taylor exapnds
$s$ in local coordinates at $x$:
\[ T^l(x)(s)=(s(x),s'(x),..,s^{(l)}(x))\in \Pr^l_{X/S}(\E)(x).\]

Assume $f:Y\rightarrow X$ is a smooth morphism of schemes over
$S$. 
\begin{proposition} \label{commutative} Let $\E$ be a finite rank locally free
  $\O_X$-module. 
There is for all $l\geq 1$ a commutative diagram of exact sequences of
locally free $\O_Y$-modules
\[
\diagram 0 \rto & \sym^l(f^*\Omega^1_{X/S}\otimes) f^*\E \rto \dto &
f^*\Pr^l_{X/S}(\E) \dto \rto & f^*\Pr^{l-1}_{X/S}(\E) \rto \dto & 0 \\
0 \rto & \sym^l(\Omega^1_{Y/S}\otimes) f^*\E \rto  &
\Pr^l_{X/S}(f^*\E) \rto & \Pr^{l-1}_{X/S}(f^*\E) \rto  & 0 \\
\enddiagram.
\]
\end{proposition}
\begin{proof}
See \cite{maa4}, Proposition 2.3.
\end{proof}

\begin{corollary} Assume $U\subseteq X$ is an open subscheme. It
  follows there is an isomorphism
\[ \Pr^l_{X/S}(\E)|_U\cong \Pr^l_{U/S}(\E|_U) \]
of $\O_U$-modules.
\end{corollary}
\begin{proof} The inclusion $i:U\rightarrow X$ is a smooth morphism
  over $S$ hence the result follows from Proposition \ref{commutative}
  and an induction.
\end{proof}

\begin{example} Jet bundles on affine schemes. \end{example}
Assume $S=\Spec(A)$ and  $U=\Spec(B)\subseteq X$ is an open affine
subscheme of $X$
and $\E|_U$ the sheafication of a locally free $B$-module $E$. It
follows $\Pr^l_{X/S}(\E)|_U$ is the sheafification of the $B$-module
\[ P^l_{B/A}(E)=B\otimes_A B/I^{l+1}\otimes_B E \]
where $I\subseteq B\otimes_A B$ is the ideal of the diagonal.

\begin{example} Jet bundles on quotients.\end{example}

Let $V$ be a finite
dimensional vector space over $K$. Let $H\subseteq G \subseteq
\GL(K,V)$ be closed subgroups.
There is an action of $H$ on $G$  defined at $K$-rational points by 
\[ \sigma: H(K)\times G(K) \rightarrow G(K) \]
\[ \sigma(h,g)=gh^{-1}.\]

The following holds: There is a \emph{quotient morphism}
\begin{align}
&\label{q1}\pi:G\rightarrow G/H 
\end{align}
and $G/H$ is a smooth quasiprojective scheme of finite type over $K$.

Let \oghmod denote the category of locally free $\O_{G/H}$-modules
of finite rank with a $G$-linearization and let \hmod denote the
category of finite dimensional $H$-modules.
There is an exact equivalence of categories
\begin{align}
&\label{q2}F: \text{\oghmod} \cong \text{\hmod}
\end{align}
defined as follows: Assume $(\E,\theta) \in$ \oghmod. It follows the
$G$-linearization $\theta$ induce an $H$-module structure 
\[ \rho(\theta): H\rightarrow \GL(\E(x))  \]
where $x\in G/H$ is the $K$-rational point defined by the identity
$e\in G$. The functor $F$ is defined by $F((\E,\theta))=(\rho(\theta),\E(x))$.
For a proof of the claims \ref{q1} and \ref{q2} see \cite{jantzen}.

In the
following when we speak of a locally free sheaf $(\E,\theta)\in$
\oghmod, and we take the fiber $\E(x)$ we will use equivalence
\ref{q2} when we refer to its $H$-module structure $\rho(\theta)$.

Let $X=G/H$ and let $Y=X\times_K X$. Let $p,q:Y\rightarrow X$ be the
canonical projection morphisms. Let $(\E,\theta)\in $\oghmod. We get from
Sequence \ref{longexact} a long exact sequence of quasi coherent
sheaves:
\begin{align}
&\label{longGH} 0\rightarrow p_*(\I_{\Delta}^{l+1}\otimes
q^*\E)\rightarrow
\H^0(G/H,\E)\otimes \O_{G/H} \rightarrow \Pr^l_{G/H}(\E) \rightarrow
\end{align}
\[\R^1p_*(\I_{\Delta}^{l+1}\otimes q^*\E)\rightarrow \R^1p_*q^*\E
\rightarrow 0 \]
Here
\[ \R^1p_*(\O_{\Delta^{l+1}}\otimes q^*\E)=0 \]
since $\O_{\Delta^{l+1}}\otimes q^*\E$ is supported on the diagonal.

There is a well defined left $G$-action on $Y$ - the
diagonal action - and since higher direct images preserve
$G$-linearizations it follows the sequence \ref{longGH} is an exact
sequence of $G$-linearized sheaves with morphisms preserving the
$G$-linearization. Since $G/H$ is a homogeneous space for the
$G$-action it follows the terms in the sequence \ref{longGH} are
locally free.

\begin{proposition} \label{glinearization} 
Assume $(\E,\theta)$ is a $G$-linearized locally free $\O_{G/H}$-module of
rank $e$. The following holds for all $l\geq 1$:
There is an exact functor
\begin{align}
&\label{pp0} \Pr^l_{G/H}:\text{\oghmod}\rightarrow \text{\oghmod}.
\end{align}
There is an exact sequence of $G$-linearized locally free $\O_{G/H}$-modules
\begin{align}
&\label{pp1} 0\rightarrow \sym^l(\Omega^1_{G/H})\otimes \E\rightarrow
\Pr^l_{G/H}(\E)\rightarrow \Pr^{l-1}_{G/H}(\E)\rightarrow 0 .
\end{align}
The Taylor morphism
\begin{align}
&\label{pp3}T^l:\H^0(G/H,\E)\otimes \O_{G/H} \rightarrow \Pr^l_{G/H}(\E)
\end{align}
preserves the $G$-linearization.
Assume $dim(G/H)=n$. The following holds:
\begin{align}
&\label{pp2} rk(\Pr^l_{G/H}(\E))=e\binom{n+l}{n}.
\end{align}
\end{proposition}
\begin{proof} 
Assume $\phi:(\E,\theta)\rightarrow (\F,\eta)$ is a morphism in
\oghmod. Let $p,q:Y\rightarrow G/H$ be the canonical projection maps
and let $Y$ have the diagonal $G$-action. It follows $q^*$ and $p_*$
preserve the $G$-linearization. We get on $Y$ a commutative diagram of
short exact sequences of morphisms of $G$-linearized sheaves
\[
\diagram 0 \rto & \I^{l+1}_\Delta \otimes q^*\E \rto \dto &
\O_Y\otimes q^*\E \rto \dto & \O_{\Delta^l}\otimes q^* \E \rto \dto &
0 \\
0 \rto & \I^{l+1}_\Delta \otimes q^*\F \rto & \O_Y \otimes q^*\F \rto
& \O_{\Delta^l}\otimes q^*\F \rto & 0
\enddiagram.
\]
Since $p_*$ preserves the $G$-linearization it follows we get a
morphism
\[ \Pr^l(\phi):\Pr^l_{G/H}(\E)\rightarrow \Pr^l_{G/H}(\F) \]
of sheaves preserving the $G$-linearization. One checks for two composable
morphisms 
\[ \phi,\psi\in \text{\oghmod} \] 
it follows 
\[\Pr^l(\phi\circ \psi)=\Pr^l(\phi)\circ \Pr^l(\psi)\]
hence the existence of the functor in claim \ref{pp0} is clear. Since $\Pr^l_{G/H}(\E)\cong
\Pr^l_{G/H}\otimes \E$ and $\Pr^l_{G/H}$ is locally free it follows
the functor is exact. It follows claim \ref{pp0} is proved.

For a proof of the existence of the sequence \ref{pp1} see \cite{maa4}, Proposition
  2.2. By the argument above it follows $\Pr^l_{G/H}(\E)$ has a
  canonical $G$-linearization since $p_*$ preserves the
  $G$-linearization. There is on $Y$ a commutative diagram of exact
  sequences of $\O_Y$-modules with a $G$-linearization
\[
\diagram 0 \rto & \I_{\Delta}^{l+1}  \rto \dto & \O_Y \dto \rto &
\O_{\Delta^l} \rto \dto & 0 \\
0 \rto & \I_{\Delta}^{l}  \rto  & \O_Y  \rto & \O_{\Delta^{l-1}} \rto  & 0
\enddiagram.
\]
When we apply the functor $p_*(-\otimes q^*\E)$ we get a
commutative diagram of maps of $G$-linearized sheaves
\[
\diagram \H^0(G/H,\E)\otimes \O_{G/H} \rto^{T^l} \dto & \Pr^l_{G/H}(\E) \dto \\
         \H^0(G/H,\E)\otimes \O_{G/H} \rto^{T^{l-1}} & \Pr^{l-1}_{G/H}(\E) 
\enddiagram
\]
hence the natural morphism $\Pr^l_{G/H}(\E)\rightarrow
\Pr^{l-1}_{G/H}(\E)$ is a morphism preserving the
$G$-linearization. It also follows the Taylor morphism preserves the
$G$-linearization.
It follows the sequence \ref{pp1} is an exact sequence of
$G$-linearized sheaves. Claim \ref{pp3} also follows from this
argument. We have proved claim \ref{pp1} and \ref{pp3}.

We prove \ref{pp2}: Assume $l=1$. We get an exact sequence
\[ 0\rightarrow \Omega^1_{G/H}\otimes \E \rightarrow
\Pr^1_{G/H}(\E)\rightarrow \E \rightarrow 0 \]
of locally free $\O_{G/H}$-modules. We know $rk(\Omega^1_{G/H})=n $ and
$rk(\E)=e$ hence
\[ rk(\Pr^1_{G/H}(\E)=ne+e=e\binom{n+1}{n}.\]
 Assume $rk(\Pr^{l-1}_{G/H}(\E)=e\binom{n+l-1}{n}$.
There is an exact sequence
\[ 0\rightarrow \sym^l(\Omega^1_{G/H})\otimes \E \rightarrow
\Pr^l_{G/H}(\E)\rightarrow \Pr^{l-1}_{G/H}(\E) \rightarrow 0 \]
hence
\[rk(\Pr^l_{G/H}(\E)=rk( \sym^l(\Omega^1_{G/H})\otimes \E
)+rk(\Pr^{l-1}_{G/H}(\E)
)=\]
\[ e\binom{n+l-1}{n-1}+e\binom{n+l-1}{n}=e\binom{n+l}{n}\]
and the Proposition is proved.
\end{proof}

Assume $(\E ,\theta)\in$ \oghmod  and consider $\Pr^l_{G/H}(\E)$. It
follows from Proposition \ref{glinearization}, claim \ref{pp1}  there is a canonical
$G$-linearization on $\Pr^l_{G/H}(\E)$. Hence we may via
equivalence \ref{q2} consider its corresponding $H$-module
$\Pr^l_{G/H}(\E)(x)$. When we speak of the $H$-module
$\Pr^l_{G/H}(\E)(x)$ we will always refer to this structure.

\begin{proposition} \label{pmod} Let $\m\subseteq \O_{G/H}$ be the
  ideal sheaf of the point $x\in G/H$. There are for all $i\geq 0$ isomorphisms of
  $H$-modules
\begin{align}
&\label{hmod1} \R^ip_*(\I_{\Delta^{l+1}}\otimes q^*\E)(x)\cong \H^i(G/H,\m^{l+1}\E) \\
&\label{hmod2} \R^ip_*q^*\E(x)\cong \H^i(G/H,\E).
\end{align}
\end{proposition}
\begin{proof} 
In the following proof when we consider a locally free sheaf
$(\E,\theta)\in$ \oghmod we will use equivalence \ref{q2} to induce
an $H$-module structure on $\E(x)$. 

Let $p:Y=G/H\times G/H\rightarrow G/H$ be defined by $p(x,y)=x$. It follows
$p^{-1}(x)\cong G/H$ and we get a fiber diagram
\[
\diagram G/H\cong p^{-1}(x)\rto^j \dto^{\pi} & G/H\times G/H \dto^p \\
               \Spec(\kappa(x)) \rto^i & G/H
\enddiagram \]
where $i(y)=(x,y)$. There is on $Y$ an exact sequence
\[ 0\rightarrow \I_\Delta^{l+1}\rightarrow \O_{Y}\rightarrow \O_{\Delta^l}\rightarrow 0\]
of $G$-linearized sheaves. Here $Y$ is equipped with the diagonal action. Let $q:Y\rightarrow G/H$ be
defined by $q(x,y)=y$. By \cite{HH}, chapter III, section 12 and
equivalence \ref{q2} we get natural maps
\[ \phi^i:\R^ip_*(\I_\Delta^{l+1}\otimes q^*\E)(x)\rightarrow
\R^i\pi_*(j^*(\I_\Delta^{l+1}\otimes q^*\E))=\H^i(G/H,\m^{l+1}\E) \]
of $H$-modules.  Let for any $\O_{Y}$-module $\F$
\[ h^i(y,\F)=dim_{\kappa(y)}\H^i(p^{-1}(y),\F_y) ,\]
where  $\F_y$ is the restriction of $\F$ to
$p^{-1}(y)$. It follows
\[ h^i(y,\I_\Delta^{l+1}\otimes q^*\E)=dim_{\kappa(y)}\H^i(p^{-1}(y),\I_\Delta^{l+1}\otimes q^*\E) \]
is a constant function in $y$ for $i\geq 0$. Hence from \cite{HH},
chapter III, Corollary 12.9 it follows $\phi^i$ is an isomorphism of
$H$-modules for all $i$.
We have proved \ref{hmod1}. Isomorphism \ref{hmod2} follows by a
similar argument and the Proposition is proved.
\end{proof}

We get from Proposition \ref{pmod}, the exact sequence \ref{longGH}
and the equivalence \ref{q2} a long exact sequence of $H$-modules
\[ 0\rightarrow \H^0(G/H,\m^{l+1}\E)\rightarrow
\H^0(G/H,\E)\rightarrow \Pr^l_{G/H}(\E)(x)\rightarrow \]
\[ \H^1(G/H,\m^{l+1}\E)\rightarrow \H^1(G/H,\E) \rightarrow 0 \]

\begin{example} Jet bundles on projective space. \end{example}

Assume in the following $V$ is a $K$-vector space of dimension $n+1$
and the characteristic of $K$ is zero. Let $L\subseteq V$ be a line  and let $P\subseteq G=\SL(V)$ be the
parabolic subgroup fixing $L$. It follows $G/P\cong \p(V^*)=\p$ - the
projective space of lines in $V$. Let $\O_\p(1)$ be the tautological
line bundle on $\p$ and let $\O_\p(d)=\O_\p(1)^{\otimes d}$. It
follows $\O_\p(d)$ has a canonical $G$-linearization and
$\Pic^G(\p)\cong \mathbf{Z}$. It follows the jet bundle
$\Pr^l_{\p}(\O_\p(d))$ has a canonical $G$-linearization for all
integers $d$. We may via equivalence \ref{q2} consider its $P$-module
$\Pr^l_{\p}(\O_\p(d))(x)$.

\begin{theorem} \label{projectivespace} There is for all $1\leq l <d$ an isomorphism
\[ \Pr^l_\p(\O_\p(d))(x)\cong \sym^l(V^*)\otimes \sym^{d-l}(L^*) \]
of $P$-modules.
\end{theorem}
\begin{proof}
The result is proved in \cite{maa1}, Theorem 2.4.
\end{proof}

Let $\pi:\p\rightarrow Y=\Spec(K)$ be the structure morphism and let
$\pi^*\sym^l(V^*)$ be the pull back of the $G$-linearized
$\O_Y$-module $\sym^l(V^*)$.

\begin{corollary} \label{split} There is for all $1\leq l < d$ an isomorphism
\[ \Pr^l_\p(\O_p(d))\cong \O_\p(d-l)\otimes \pi^*\sym^l(V^*) \]
of $\O_\p$-modules with a $G$-linearization.
\end{corollary}
\begin{proof} Using equivalence \ref{q2} it follows the $P$-module
  corresponding
to $\O_\p(d-l)$ is $\sym^{d-l}(L^*)$. It follows from Theorem \ref{projectivespace}
$\Pr^l_\p(\O_p(d))$ and $\O_\p(d-l)\otimes \pi^*\sym^l(V^*)$
have isomorphic $P$-modules. The Corollary follows since \ref{q2} is
an equivalence of categories.
\end{proof}
We get for all $1\leq l <d$ get the following formula
\begin{align}
&\label{splittingtype} \Pr^l_\p(\O_\p(d))\cong \oplus^{\binom{n+l}{n}}\O_\p(d-l)
\end{align}
expressing the splitting type of the jet bundle as abstract locally
free $\O_\p$-module. This follows from Corollary \ref{split} since the $\O_\p$-module $\pi^*\sym^l(V^*)$
corresponds to the trivial rank $\binom{n+l}{n}$ $\O_\p$-module.
The formula \ref{splittingtype} is well known (see \cite{maa1}, \cite{maa2}, \cite{maa4},
\cite{perk},\cite{pien}, \cite{diro} and \cite{somm} for other proofs.)

\section{On surjectivity of the Taylor morphism}

In this section we study the Taylor morphism for a class of invertible
sheaves on the grassmannian. We prove the Taylor morphism is
surjective in most cases. 
 
The strategy of the proof is as follows: First prove the Taylor morphism is
surjective for most invertible sheaves on projective space. 
Embed the grassmannian into projective space using the Pl\"{u}cker embedding.
Any invertible sheaf on the grassmannian is the pull back of an
invertible sheaf on projective space via the Pl\"{u}cker embedding.
Using the fact the grassmannian is projectively normal in the Pl\"{u}cker embedding
and general properties of jet bundles we prove the Taylor morphism is
surjective for most invertible sheaves on the grassmannian.

\textbf{Notation:} Let $K$ be any field and let $W\subseteq V$ be
$K$ vector spaces of dimensions $m$ and $m+n$. Let $G=\SL(V)$ and let
$P\subseteq G$ be the subgroup fixing $W$. Let $X=G/P=\gr(m,m+n)$. Let
\[ i:X\rightarrow \p(\wedge^m V^*) \]
be the Pl\"{u}cker embedding and let $\O_X(d)=i^*\O(1)^{\otimes
  d}$. The grassmannian $\gr(m,m+n)$ has dimension $mn$.
There is a bijection
\begin{align} \{\text{$K$-rational points $x\in G/P$}\} \cong
  \{\text{$m$-planes $W_x\subseteq V$}\}
\end{align}
of sets.
Let \pmod  be the category of finite dimensional $P$-modules and morphisms
and let \ogpmod be the category of $G$-linearized locally free $\O_{G/P}$-modules
of finite rank and morphisms.
Recall from the previous section there is an exact equivalence of categories
\begin{align} \label{equivalence} \text{\pmod} \cong \text{\ogpmod}.
\end{align}

Let $Y=X\times_K X$ and ler $p,q:Y\rightarrow X$ be the
canonical projection morphisms. Let $\Delta(X)\subseteq Y$ be the
diagonal embedding of $X$ and let $\I_\Delta \subseteq \O_Y$ be the
ideal of $\Delta(X)$.
By Definition \ref{jets} it follows
\[ \Pr^l_X(\O_X(d))=p_*(\O_Y/\I_\Delta^{l+1}\otimes q^*\O_X(d)) \]
is the \emph{$l$'th sheaf of jets } of $\O_X(d)$.
It follows from Proposition \ref{glinearization}, claim \ref{pp3}
\[ rk(\Pr^l_X(\O_X(d))) =\binom{mn+l}{mn}.\]
When it is clear from the context we will write $\Pr^l(\O(d))$ instead
of $\Pr^l_X(\O_X(d))$.

The invertible sheaf $\O_{X}(d)$ has a unique $G$-linearization for
all $d\in \mathbf{Z}$. Let the product $Y=G/P\times_K G/P$ have
the diagonal $G$-action. 

Recall there is an exact sequence of $G$-linearized
sheaves of $\O_Y$-modules
\[ 0 \rightarrow \I_\Delta^{l+1}\rightarrow \O_Y \rightarrow \O_{\Delta^l}
\rightarrow 0 \]
on $Y$. The functor $p_*(-\otimes q^*\O_X(d)) $ is left exact and
preserves the $G$-linearization hence we get when we use the formalism
of derived functors a long exact sequence of locally free $G$-linearized sheaves
\begin{align} \label{ex1}0\rightarrow p_*(\I_\Delta^{l+1}\otimes q^*\O_X(d))\rightarrow p_*q^*\O_X(d)
\rightarrow \Pr^l_X(\O_X(d)) \rightarrow \end{align}
\[ \R^1p_*(\I_\Delta^{l+1}\otimes q^*\O_X(d))\rightarrow
\R^1p_*q^*\O_X(d)\rightarrow  0 .\]
Here  $\R^1p_*(\O_{\Delta^l}\otimes q^*\O_X(d))=0$ since
$\O_{\Delta^l}\otimes q^*\O_X(d)$ is a sheaf supported on the diagonal.

Let $x\in G/P$ be the $K$-rational point defined by the class of the
identity element $e\in G$.  When we take the fiber at $x$ of the sequence \ref{ex1} and apply
Proposition \ref{pmod} we get the following exact sequence of finite
dimensional $P$-modules
\begin{align} \label{fiber1}
0\rightarrow \H^0(X,\m^{l+1}\O_X(d))\rightarrow \H^0(X,\O_X(d))\rightarrow
\Pr^l_X(\O_X(d))(x)\rightarrow 
\end{align}
\[ \H^1(X,\m^{l+1}\O_X(d))\rightarrow \H^1(X,\O_X(d)) \rightarrow \cdots
\]
By Kodaira's Vanishing Theorem it follows $\H^1(X,\O_X(d))=0$ when $d\geq
1$. It follows we get an exact sequence
\begin{align} \label{fiber2}
0\rightarrow \H^0(X,\m^{l+1}\O_X(d))\rightarrow \H^0(X,\O_X(d))\rightarrow^{T^l}
\Pr^l_X(\O_X(d))(x)\rightarrow
\end{align}
\[ \H^1(X,\m^{l+1}\O_X(d))\rightarrow 0 \]
of finite dimensional $P$-modules.
Since \ref{pmod} is an equivalence of categories, we get an exact
sequence of locally free $G$-linearized sheaves
\begin{align} \label{ex2}0\rightarrow p_*(\I\Delta^{l+1}\otimes
  q^*\O_X(d))\rightarrow \H^0(X,\O_X(d))\otimes \O_X 
\rightarrow \Pr^l_X(\O_X(d)) \rightarrow \end{align}
\[ \R^1p_*(\I_\Delta^{l+1}\otimes q^*\O_X(d))\rightarrow 0 .\]

\begin{example} Taylor maps and invertible sheaves on projective space.\end{example}

Let $E$ be an $n+1$-dimensional $K$-vector space and let
$\p=\p(E^*)$ be the projective space of lines in $E$. Let $\O_\p(1)$
be the tautological line bundle on $\p$ and let $\O_\p(d)=\O_\p(1)^{\otimes d}$. 

\begin{lemma} \label{projective} The Taylor morphism
\[ T^l:\H^0(\p,\O_\p(d))\otimes \O_\p \rightarrow \Pr^l_\p(\O_\p(d)) \]
is surjective for all $1\leq l \leq d$.
\end{lemma}
\begin{proof} Let $E=K\{e_0,..,e_n\}$ and $E^*=K\{x_0,..,x_n\}$. It
  follows
\[ \p(E^*)=\Proj(\sym_K(E^*))=\Proj(K[x_0,..,x_n]).\]
Let $U_0=D(x_0)=\Spec(K[t_1,..,t_n])$ where $t_i=x_i/x_0$. There is an
isomorphism 
\[ \p(E^*)\cong \SL(E)/P \]
where $P$ is the subgroup of elements fixing a line in $E$. 
Because the Taylor morphism is a map of $\SL(E)$-linearized sheaves we
may check surjectivity by restricting to the fiber of $T^l$ at $x$. We
restrict $T^l$ to the open set $U_0$:
\[ T^l|_{U_0}:K[t_1,..,t_n]\otimes \H^0(\p,\O_\p(d))\rightarrow
\Pr^l_{U_0}(\O_\p(d)|_{U_0}).\]
We get a map
\[ T^l|_{U_0}:K[t_i]\otimes \H^0(\p,\O_\p(d))\rightarrow K[t_i]\otimes 
\{dt_1^{d_1}\cdots dt_n^{d_n}\otimes x_0^d : 0\leq \sum d_i \leq n\} \]
of left $K[t_i]$-modules.
Assume
\[ s= x_0^{d_0}x_1^{d_1}\cdots x_n^{d_n} \]
with $\sum d_i=d$ is a global section of $\O_\p(d)$. It follows
$d_0=d-d_1-\cdots -d_n$ On $U_0$ we may
write
\[ s=x_0^{d-d_1-\cdots -d_n}x_1^{d_1}\cdots x_n^{d_n}=t_1^{d_1}\cdots
t_n^{d_n}x_0^d.\]
By definition
\[ T^l(s)=1\otimes t_1^{d_1}\cdots t_n^{d_n}x_0^d=\]
\[ (t_1+dt_1)^{d_1}\cdots (t_n+dt_n)^{d_n}\otimes x_0^d .\]
The point $x$ is defined by $t_1=\cdots =t_n=0$ hence when we restrict
$T^l$ to the fiber at $x$ we get the map
\[ T^l(x):\H^0(\p,\O_\p(d))\rightarrow \Pr^l_\p(\O_\p(d))(x) \]
defined by
\[ T^l(x)(s)=dt_1^{d_1}\cdots dt_n^{d_n}\otimes x_0^d.\]
Assume $\omega=dt_1^{d_1}\cdots dt_n^{d_n}\otimes x_0^d \in \Pr^l_\p(\O_\p(d))(x)$
with $0\leq \sum d_i \leq n$. It follows $d-\sum d_i \geq d-n \geq 0$.
Let $d_0=d-\sum d_i$. It follows $d_0\geq 0$ and $d_0+\sum d_i=d$. It
follows $s=x_0^{d_0}x_1^{d_1}\cdots x_n^{d_n}\in \H^0(\p,\O_\p(d))$ and
\[ T^l(x)(s)=\omega \]
and the Proposition is proved.
\end{proof}

Recall we get on projective space $\p$ an exact sequence of
$\SL(E)$-linearized locally free sheaves
\begin{align} \label{proj} 0\rightarrow p_*(\I_\Delta^{l+1}\otimes q^*\O_\p(d))\rightarrow \H^0(\p,
\O_\p(d))\otimes \O_\p \rightarrow^{T^l} \Pr^l_\p(\O_\p(d))
\rightarrow \end{align}
\[ \R^1p_*(\I_\Delta^{l+1}\otimes q^*\O_\p(d)) \rightarrow  \R^1p_*(q^*\O_\p(d)) \rightarrow 0\]
when $1\leq l \leq d$.

\begin{corollary}  There is an equality
\[ \R^1p_*(\I_\Delta^{l+1}\otimes q^*\O_\p(d))=0 \]
when $1\leq l \leq d$.
\end{corollary}
\begin{proof} Sequence \ref{proj} remain exact when we take the fiber
  at $x\in G/P$. Via Kodaira's Vanishing Theorem and Proposition \ref{pmod}, claim \ref{hmod2}
  the final term becomes
\[ \R^1p_*(q^*\O_\p(d))(x)=\H^1(\p, \O_\p(d))=0 \]
when $d\geq 1$. It follows $\R^1p_*q^*\O_p(d)=0$.
Since $T^l$ is surjective when $1\leq l \leq d$ the
Corollary follows.
\end{proof}

We get on $\p=\p(E^*)$ an exact sequence of $\SL(E)$-linearized locally free sheaves
\begin{align} \label{proj0} 0\rightarrow p_*(\I_\Delta^{l+1}\otimes
  q^*\O_\p(d))\rightarrow \H^0(\p,
\O_\p(d))\otimes \O_\p \rightarrow^{T^l} \Pr^l_\p(\O_\p(d))
\rightarrow 0 \end{align}
when $1\leq l \leq d$.

\begin{example} Surjectivity of Taylor maps for projectively normal schemes. \end{example}

\begin{lemma} \label{surjective} Let $i:Z\rightarrow W$ be a closed immersion of
  schemes and let $\E$ be a locally free $\O_W$-module. There is a
  canonical surjection
\[ \phi:i^*\Pr^l_W(\E)\rightarrow \Pr^l_Z(i^*\E) \]
of $\O_Z$-modules.
\end{lemma}
\begin{proof} Assume $Z=\Spec(A/\a)$, $W=\Spec(A)$ and $\E=\tilde{E}$
  where $E$ is a locally free $A$-module.
Let $\Pr^l_W(\E)$ be the sheafification of $A\otimes A/J^{l+1}\otimes
E$ where $J\subseteq A\otimes A$ is the ideal of the diagonal.
Moreover $\Pr^l_Z(i^*\E)$ the sheafification of $(A/\a)\otimes
(A/\a)/\tilde{J}^{l+1}\otimes (E/IE)$
where $\tilde{J}\subseteq A/\a \otimes A/\a$ is the ideal of the diagonal.
There is an isomorphism between $i^*\Pr^l_W(\E)$ and the
sheafification of
\[A\otimes A/J^{l+1}\otimes (E/\a E) .\]
In this case the map $\phi$ is the sheafification of the canonical map
\[ f:A\otimes A/J^{l+1}\otimes (E/\a E) \rightarrow (A/\a )\otimes
(A/\a )/\tilde{J}^{l+1}\otimes (E/\a E) \]
defined by
\[ f(x\otimes y \otimes \overline{e})=\overline{x}\otimes
\overline{y}\otimes \overline{e} .\]
It follows $\phi$ is a surjective map of sheaves. This construction
glue to give a morphism for any closed immersion and the Lemma is
proved.
\end{proof}

Assume $i:Z\rightarrow \p^N_K$ is an embedding of a projective scheme $Z$
and assume $Z$ is projectively normal in $i$. Assume $\F$ is a coherent $\O_\p$-module with 
\[ T^l:\H^0(\p,\F)\otimes \O_\p \rightarrow \Pr^l_\p(\F) \]
surjective and let $\E=i^*\F$.

\begin{theorem} \label{mainsurjective} The Taylor morphism
\[ T^l:\H^0(Z,\E)\otimes \O_Z \rightarrow \Pr^l_Z(\E) \]
is surjective.
\end{theorem}
\begin{proof} 
Pull $T^l$ back to $Z$ via $i$
to get a surjective morphism of sheaves
\[ i^*(T^l):\H^0(\p,\F)\otimes \O_Z\rightarrow i^*\Pr^l_\p(\F)
.\]
We get a commutative diagram of maps of sheaves
\[
\diagram \H^0(\p,\F) \otimes \O_Z \rto^{i^*(T^l)} \dto^u &
i^*\Pr^l_\p(\F) \dto^{\phi} \\
\H^0(Z,\E)\otimes \O_Z \rto^{T^l} & \Pr^l_Z(\E)
\enddiagram .\]
The map $u$ is surjective since the $Z$ is projectively
normal in the embedding $i$ and $\phi$ is surjective by Lemma
\ref{surjective} and the Theorem is proved.
\end{proof}

\begin{example} Higher cohomology groups of coherent sheaves on $G/P$.
\end{example}

Assume $Z=G/P$ where $P\subseteq G$ is a parabolic subgroup of a semi simple
algebraic group. Assume $i:G/P\rightarrow \p^N_K$ is an embedding
where $G/P$ is projectively normal in $i$.
Let $\O_{G/P}(d)=i^*\O(d)$.

\begin{proposition} \label{nadel} Assume $1\leq l <d$. If $\R^1p_*q^*\O_Z(d)=0$ it follows
\[\R^1p_*(\I^{l+1}\otimes q^*\O_Z(d))=0.\]
\end{proposition}
\begin{proof} Let $\I\subseteq \O_{G/P\times G/P}$ be the ideal sheaf of
  the diagonal. We get an exact sequence of sheaves on $G/P$:
\[ 0\rightarrow p_*(\I^{l+1}\otimes q^*\O_{G/P}(d))\rightarrow
p_*q^*\O_{G/P}(d)\rightarrow^{T^l} \Pr^l(\O_{G/P}(d)) \rightarrow^\phi \]
\[ \R^1p_*(\I^{l+1}\otimes q^*\O_{G/P}(d)) \rightarrow^\psi
\R^1p_*q^*\O_{G/P}(d)\rightarrow \cdots \]
Assume $\R^1p_*q^*\O_Z(d)=0$. Since $G/P$ is projectively normal in
$i$ it follows from Theorem \ref{mainsurjective} the Taylor morphism
$T^l$ is surjective. 
It follows
\[ Im(T^l)=\Pr^l(\O_Z(d))=Ker(\phi).\]
It follows $\phi(\Pr^l(\O_Z(d))=0$ hence $Im(\phi)=0=Ker(\psi)$. 
We get an injection
\[ \R^1p_*(\I^{l+1}\otimes q^*\O_Z(d))\rightarrow \R^1p_*q^*\O_Z(d)=0
\]
and the Proposition is proved.
\end{proof}

Assume the hypothesis from Proposition \ref{nadel} is satisfied.

\begin{corollary} \label{nadelvanishing} There is for $1\leq l <d$ an equality
\[ \H^1(G/P,\m^{l+1}\O_Z(d))=0.\]
\end{corollary}
\begin{proof} By Proposition \ref{nadel} it follows
\[ \R^1p_*(\I^{l+1}\otimes q^*\O_Z(d))=0.\]
Take the fiber at the distinguished point $x\in G/P$ to get the
equality
\[ \R^1p_*(\I^{l+1}\otimes
q^*\O_Z(d))(x)=\H^1(G/P,\m^{l+1}\O_Z(d))=0\]
and the Corollary is proved.
\end{proof}

Note: Corollary \ref{nadelvanishing} is a special case of Nadel's Vanishing
Theorem on the vanishing of higher cohomology groups of coherent
sheaves on flag varieties. Proposition \ref{nadel} shows there is an
equation
\[ \R^1p_*(\I^{l+1}\otimes q^*\O_Z(d))=0.\]
The coherent sheaf $\R^1p_*(\I^{l+1}\otimes q^*\O_Z(d))$ is a locally
free $\O_{G/P}$-module and its fiber at $x$ is the cohomology group
\[ \H^1(G/P,\m^{l+1}\O_Z(d)).\]
Hence Corollary \ref{nadelvanishing} gives a geometric proof of the
vanishing of $\H^1$ for a class of coherent sheaves on flagvarieties
over an arbitrary field $K$ of characteristic zero.

Note: By examining the proof we see the result is true over an
arbitrary field.

\begin{example} Taylor maps for line bundles on the grassmannian.\end{example}

Let $X=\SL(V)/P=\gr(m,m+n)$ be the grassmannian of $m$-planes in $V$
and let $\O_X(d)$ be the line bundle coming from the Pl\"{u}cker
embedding.

\begin{corollary} \label{surjective} The Taylor morphism
\[ T^l:\H^0(X,\O_X(d))\otimes \O_X \rightarrow \Pr^l_X(\O_X(d)) \]
is surjective when $1\leq l \leq d$.
\end{corollary}
\begin{proof} The Corollary follows from Theorem \ref{mainsurjective}
since the grassmannian is projectively normal in
the Plucker embedding. 
\end{proof}

On $X=\gr(m,m+n)$ we get an exact sequence of $\SL(V)$-linearized
sheaves

\begin{align} \label{proj1} 0\rightarrow p_*(\I_\Delta^{l+1}\otimes
  q^*\O_X(d))\rightarrow \H^0(X, \O_X(d))\otimes \O_X \rightarrow^{T^l} \Pr^l_X(\O_X(d))
\rightarrow \end{align}
\[ \R^1p_*(\I_\Delta^{l+1}\otimes q^*\O_X(d)) \rightarrow
\R^1p_*(q^*\O_X(d)) \rightarrow 0\]
when $1\leq l \leq d$.

\begin{corollary} \label{vanishinggrass} On $X=\gr(m,m+n)$ there is an equality
\[ \R^1p_*(\I_\Delta^{l+1}\otimes q^*\O_X(d))=0\]
when $1\leq l \leq d$.
\end{corollary}
\begin{proof} The sequence \ref{proj1} remain by the equivalence
  \ref{q2} exact when we take the fiber at $x$. We get by Proposition
  \ref{pmod}, claim \ref{hmod2}
\[ \R^1p_*(q^*\O_X(d))(x)=\H^1(X,\O_X(d))=0 \]
when $d\geq 1$ by Kodaira's Vanishing Theorem. It follows $\R^1p_*q^*\O_X(d)=0$.
By Theorem \ref{surjective}
the Taylor map $T^l$ is surjective and the Corollary follows since
\ref{proj1} is an exact sequence.
\end{proof}

%\begin{example} \label{nadel}  Nadel's Vanishing Theorem.
%\end{example}
%When we take the fiber of $\R^1p_*(\I^{l+1}\otimes q^*\O_X(d))$ at $x$
%we get from Corollary \ref{vanishinggrass} and Proposition \ref{pmod} an equality
%\[ \H^1(X,\m^{l+1}\O_X(d))=0 \]
%in the case when $1\leq l \leq d$. This is a special case of 
%Nadel's Vanishing Theorem on the vanishing of the cohomology of a
%class of coherent sheaves on the grassmannian. 
%It may in the complex analytic case be proved using the theory of multiplier ideal sheaves.

We get an exact sequence of $\SL(V)$-linearized sheaves
\begin{align} \label{proj2} 0\rightarrow p_*(\I_\Delta^{l+1}\otimes
  q^*\O_X(d))\rightarrow \H^0(X, \O_X(d))\otimes \O_X
  \rightarrow^{T^l} \Pr^l_X(\O_X(d))
\rightarrow 0 \end{align}
when $1\leq l \leq d$.

\begin{corollary} \label{pmod1} On $X=\gr(m,m+n)$ there is an exact sequence of
$P$-modules
\[ 0\rightarrow \H^0(X,\m^{l+1}\O_X(d))\rightarrow
\H^0(X,\O_X(d))\rightarrow \Pr^l_X(\O_X(d))(x)\rightarrow 0 \]
for all $1\leq l \leq d$.
\end{corollary}
\begin{proof} If we take the fiber of sequence \ref{proj2} we get via
  equivalence \ref{q2} and Proposition \ref{pmod} an exact sequence
\[ 0\rightarrow \H^0(X,\m^{l+1}\O_X(d))\rightarrow
\H^0(X,\O_X(d))\rightarrow \Pr^l_X(\O_X(d))(x)\rightarrow 0 \]
of $P$-modules and the Corollary follows.
\end{proof}

\section{On generalized Verma modules and canonical filtrations}

In this section we study the canonical filtration $\U_l(\lg)v$ associated to the
irreducible $\SL(V)$-module of global sections of a line bundle $\L\in
\Pic^G(G/P)$ on the grassmannian $G/P=\gr(m,m+n)$ of $m$-planes in an
$m+n$-dimensional vector space. We construct a basis and compute the
dimension of each term in the canonical filtration.

The strategy of the proof is as follows: There is for all $l\geq 1$ an exact sequence
\[ 0 \rightarrow ann_l(v)\otimes L_v \rightarrow \U_l(\lg)\otimes L_v
\rightarrow \U_l(\lg)v \rightarrow 0 \]
where $ann_l(v)$ is the $l$'th piece of the canonical filtration of
the annihilator ideal of $v$ - the highest weight vector of $\H^0(G/P,
\O_{G/P}(d))^*$. Here $\U_l(\lg)$ is the $l$'th piece of the canonical
filtration of the universal enveloping algebra of $\lg=Lie(G)$.
Using the theory of highest weigts and the Poincare-Birkhoff-Witt
Theorem we prove there is a vector space decomposition
\begin{align}
&\label{dec}  \U_l(\lg)=\U_l(\ln)\oplus ann_l(v) 
\end{align}
in the case when $1\leq l <d$. Here $\ln \subseteq \lg$ is a sub Lie
algebra. We then use the decomposition \ref{dec} to give a basis for
$\U_l(\lg)v$ and to calculate its dimension as vector space.

\textbf{Notation:} Let $W\subseteq V$ be vector spaces of dimension
$m$ and $m+n$. Let $e_1,..,e_{m+n}$ be a basis for $V$ and
$e_1,..,e_m$ a basis for $W$. Let $G=\SL(V)$ and $P\subseteq G$ the
parabolic subgroup fixing $W$. It follows there is a quotient morphism
\[ \pi:G\rightarrow G/P \]
and there is a canonical isomorphism
\[ G/P\cong \gr(m,m+n) \]
where $\gr(m,m+n)$ is the grassmannian of $m$-planes in $V$.
Let 
\[ i:G/P\rightarrow \p(\wedge^m V^*) \]
be the Pl\"{u}cker embedding and let $\O_X(d)=i^*\O(1)^{\otimes d}$.
Let $\lg=Lie(G)$ and $\lp=Lie(P)$. Let $\U(\lg)$ be the universal
enveloping algebra of $\lg$ and let $\U_l(\lg)$ be its canonical
filtration. Let 
\[ \Vl=\H^0(G/P,\O_{G/P}(d))^* \]
be the irreducible $\SL(V)$-module of dual global sections of
$\O_{G/P}(d)$ and let $v\in \Vl$ be the (unique up to scalars) highest
weight vector of $\Vl$.

In the following we use the notation from \cite{fulton}.
Let $\lg=\lg_{-}\oplus \lh \oplus \lg_{+}$ be the
\emph{triangular decomposition} of $\lg$ defined as follows: Elements in $\lg$ are matrices
$A$ of dimension $m+n$ with trace zero. Let $\lg_-$ be the set of
lower triangular matrices in $\lg$, $\lg_+$ the set of upper
triangular matrices in $\lg$ and $\lh$ the set of diagonal
matrices. Hence $\lh$ is the Lie algebra of diagonal matrices $A$ of
dimension $m+n$ with trace zero. It follows $\lh$ consists of matrices of the
type
\[
A=\begin{pmatrix} a_1 & 0    &  0 & \cdots & 0 & 0 \\
                   0  & a_2  &  0 & \cdots & 0 & 0 \\
                   \vdots &  \vdots & \cdots & \cdots & \vdots & \vdots \\
                  0 & 0 & 0 & \cdots & 0 & a_{m+n}  
\end{pmatrix}
\]
with $tr(A)=\sum a_i=0$. Let $\lh^*$ be the dual of $\lh$. It follows
\[ \lh^*=K\{L_1,..,L_{m+n}\}/L_1+\cdots +L_{m+n} \]
where
\[ L_i(A)=a_i.\]

Note. Given a semi simple Lie algebra $\lg$ there exist many different
Cartan decompositions. They are all conjugate under automorphisms of $\lg$.

There is since $K$ has characteristic zero an embedding of
$G$-modules
\[ \Vl\subseteq \sym^d(\wedge^mV^*)^*\cong \sym^d(\wedge^m V) .\]
Let 
\[ l^d=\sym^d(\wedge^mW)\subseteq \sym^d(\wedge^m V).\]
It follows $l^d$ is a $P$-stable vector. Let $L_{l^d}\subseteq
\sym^d(\wedge^m V)$ be the line spanned by $l^d$ and let $L_v\subseteq
\Vl$ be the line spanned by $v$.

\begin{lemma} \label{stable}There is an equality $L_v=L_{l^d}$. Moreover $v$ is the
 unique  highest weight vector for $\Vl=\H^0(X,\O_X(d))^*$ with highest weight
\[ \lambda=d(L_1+\cdots +L_m).\]
\end{lemma}
\begin{proof} By the Borel-Weil-Bott Theorem it follows $\Vl$ is an
  irreducible $G$-module. One checks there is an equality
  $L_v=L_{l^d}$ of lines, and that $v$ is a highest weight vector for $\Vl$. 
One also checks $v$ has the given weight and the Lemma follows.
\end{proof}

The line $L_v \subseteq \Vl$ is in fact the unique $P$-stable line of
$\Vl$. 

The $l$'th piece $\U_l(\lg)$ of the canonical filtration of the
enveloping algebra is a $G$-module via the adjoint representation. 
It follows $\U_l(\lg)$ is a $P$-module.
\begin{definition} Let $ann(v)\subseteq \U(\lg)$ be the left
  \emph{annihilator ideal} of the vector $v\in \Vl$:
\[ ann(v)=\{x\in \U(\lg) : x(v)=0\}.\]
 Let
  $ann_l(v)=ann(v)\cap \U_l(\lg)$ be its \emph{canonical filtration}.
\end{definition}
It follows $ann(v)$ is a left ideal in the associative ring $\U(\lg)$.

There is an exact sequence
\[ 0\rightarrow ann(v)\otimes_K L_v \rightarrow \U(\lg)\otimes_K L_v
\rightarrow \Vl \rightarrow 0 \]
of $G$-modules
and an exact sequence
\[ 0\rightarrow ann_l(v)\otimes_K L_v \rightarrow \U_l(\lg)\otimes_K L_v
\rightarrow \U_l(\lg)v \rightarrow 0 \]
of $P$-modules. Here 
\[ \U_l(\lg)v=\{ x(v): x\in \U_l(\lg) \} \subseteq \Vl \]
is the $P$-module
spanned by elements of $\U_l(\lg)$ and the vector $v$.
The $G$-module $\U(\lg)\otimes L_v$ is the \emph{generalized Verma
  module} associated to the $P$-module $L_v$. The $G$-module
$\U(\lg)\otimes_K L_v$ has a canonical filtration of $P$-modules given by 
\[ \U_l(\lg)\otimes L_v \subseteq \U(\lg)\otimes L_v .\]

\begin{definition} Let $\{\U_l(\lg)v\}_{l\geq 0}$ be the \emph{canonical
  filtration} of $\Vl$.
\end{definition}
The $P$-module $\U_l(\lg)v$ depends on the $P$-stable line $L_v$
defined by $v\in \Vl$ 
which is canonically determined by the highest weight vector $v\in
\Vl$.
It follows we get a canonical filtration
\[ \U_1(\lg)v\subseteq \cdots \subseteq \U_l(\lg)v\subseteq \Vl \] 
of $\Vl$ by $P$-modules.

The Lie algebra $\lp$ is the sub Lie algebra of $\lg$
consisting of traceless matrices $M$ with coefficients in $K$ on the following form:
\[M=
\begin{pmatrix}  A & X \\
                 0 & B \end{pmatrix} 
\]
where $A$ is an $m\times m$-matrix, $X$ is an
$n\times m$-matrix and $B$ is an $n\times n$-matrix 
such that $tr(A)+tr(B)=0$. 
Let $\ln \subseteq \lg$ be the sub Lie algebra consisting of matrices
$M$ on the following form:
\[M=
\begin{pmatrix}  0 & 0 \\
                 Y & 0 \end{pmatrix}
\]
where $Y$ is an $m\times n$-matrix with coefficients in $K$. It
follows there is a direct sum decomposition $\lg=\ln\oplus \lp$ of
vector spaces. 

\begin{proposition} The natural map
\[ f:\U(\ln)\otimes_K \U(\lp)\rightarrow \U(\lg) \]
defined by
\[ f(v\otimes w)=vw \]
is an isomorphism of vector spaces.
\end{proposition}
\begin{proof} For a proof, see \cite{dixmier}, Proposition 2.2.9.
\end{proof}

\begin{definition} Let $l\geq 1$ be an integer and define
\[ \U_l(\ln,\lp)=\sum_{i+j=l}\U_i(\ln)\otimes_K \U_j(\lp)\subseteq
\U(\ln)\otimes_K \U(\lp) .\]
\end{definition}

Assume in the following $\ln$ has $\{x_1,..,,x_E\}$ as a basis.
The line $v=l^d$ define a $\lp$-module
\[ \rho:\lp\rightarrow \End(v) \]
hence we get a map
\[ \rho:\lp \rightarrow K .\]
If $z\in \lp$ is an element defined by the matrix $M$ with 
\[M=
\begin{pmatrix}  A & X \\
                 0 & B \end{pmatrix}
\]
it follows $\rho(z)=dtr(A)$.
It follows we may choose a direct sum decomposition $\lp=\lp_v\oplus
(x)$ where $\lp_v=Ker(\rho)$ and $x\in \lp$ has the property
$\rho(x)=d$.
It follows we may choose a basis for $\lp$ on the form
$\{y_1,..,y_D,x\}$ where $\rho(y_i)=0$.

\begin{lemma} There is for all $l\geq 1$ an equality
\[ f(\U_l(\ln,\lp))=\U_l(\lg).\]
\end{lemma}
\begin{proof} It is clear $f(\U_l(\ln,\lp))\subseteq
  \U_l(\lg)$. Assume
\[ \omega=x_1^{v_1}\cdots x_E^{v_E}x^qy_1^{u_1}\cdots y_D^{u_D}\in
\U_l(\lg)\]
with $\sum_i v_i+q+\sum_j u_j\leq l$.
It follows
\[ z= x_1^{v_1}\cdots x_E^{v_E}\otimes x^qy_1^{u_1}\cdots y_D^{u_D}\in
\U_l(\ln,\lp) \]
and $f(z)=\omega$. The Lemma is proved.
\end{proof}

Let $\onep \in \U(\lp)$ be the multiplicative identity element. 

\begin{definition} Let $l\geq 1$ be an integer. Define the following:
\[ \U_l(\ln)\otimes \onep=\{ x\otimes \onep: x\in
\U_l(\ln)\} \subseteq \U_l(\ln,\lp)\]
\[ W_l=\{x\otimes w(y-\rho(y)\onep): x\in \U_i(\ln), y\in \lp,
w(y-\rho(y)\onep)\in \U_j(\lp), i+j=l\}\subseteq \U_l(\ln,\lp).\]
\end{definition}

\begin{proposition} \label{injective} The natural map
\[ \phi:\U_l(\ln)\otimes \onep \oplus W_l \rightarrow \U_l(\ln,\lp) \]
defined by
\[ \phi(x,y)=x+y \]
is an isomorphism of vector spaces for all $l\geq 1$.
\end{proposition}
\begin{proof} 
We first prove the following:
\[ \U_l(\ln)\otimes \onep \cap W_l=\{0\} \]
for all $l\geq 1$.
Assume $z\in W_l$ is a monomial with $l\geq 1$. It follows
\[ z=x\otimes w(y-\rho(y)\onep) \]
with $y\in \lp$. One of the following assertions is true:
\begin{align}
&\label{a1} z=x\otimes w y_i  \\
& \label{a2} z=x\otimes w(x-d\onep) \\
& \label{a3} z=0
\end{align}
Since $x\otimes w y_i$ and $x\otimes w(x-d\onep)$ are not in
$\U_l(\ln)\otimes \onep$ it follows $z=0$ and the claim is proved.
It follows the map $\phi$ is an injection.

We next prove the map $\phi$ is a surjection for all $l\geq 1$:
By definition
\[ \U_l(\ln,\lp)=\sum_{i+j=l}\U_{l-i}(\ln)\otimes_K \U_i(\lp) .\]
Pick a monomial
\[z=x_1^{v_1}\cdots x_E^{v_E}\otimes x^qy_1^{u_1}\cdots y_D^{u_D}\in \U_l(\ln,\lp) .\]
If $q=u_1=\cdots =u_D=0$ it follows $z\in \U_l(\ln)\otimes
\onep$. Assume
$q+\sum u_i \geq 1$ and $1\leq j \leq D$ the largest integer with
$u_j\geq 1$. It follows
\[ z=x_1^{v_1}\cdots x_E^{v_E}\otimes x^qy_1^{u_1}\cdots y_j^{u_j}=\]
\[ x_1^{v_1}\cdots x_E^{v_E}\otimes x^qy_1^{u_1}\cdots y_j^{u_j-1}(y_j-\rho(y_j)\onep)\in W_l.\]
Assume $u_1=\cdots =u_D=0$ and $q\geq 1$. We may write
\[ x^q=d^q\onep +y(x-d\onep) \]
for some $y$. We get
\[ z= x_1^{v_1}\cdots x_E^{v_E}\otimes x^q= x_1^{v_1}\cdots
x_E^{v_E}\otimes (d^q\onep +y(x-d\onep))=\]
\[x_1^{v_1}\cdots x_E^{v_E}d^q\otimes \onep+x_1^{v_1}\cdots
x_E^{v_E}\otimes y(x-d\onep)\in \U_l(\ln)\otimes \onep +W_l \]
and the Proposition is proved.
\end{proof}

Let $\oneg\in \U(\lg)$ be the multiplicative identity element.
\begin{definition} Let 
\[ char(\rho)= \{x(y-\rho(y)\oneg):y\in \lp, x\in \U(\lg) \} \] 
be the left \emph{character ideal} of $\rho$. Let 
\[ char_l(\rho)=char(\rho)\cap \U_l(\lg) \]
be the canonical filtration of $char(\rho)$.
\end{definition}
It follows $char(\rho)\subseteq \U(\lg)$ is a left ideal in the associative ring $\U(\lg)$.
By definition it follows there is an inclusion of ideals
$char(\rho)\subseteq ann(v)$. This inclusion is strict.
The inclusion $\ln \subseteq \lg$ induce a canonical injection
$\U(\ln)\subseteq \U(\lg) $ of associative rings. We get a canonical inclusion
of filtrations $\U_l(\ln)\subseteq \U_l(\lg)$ for all $l\geq 1$.

\begin{lemma} \label{inclusion} The following holds for all $l\geq $:
\begin{align}
&\label{m1}f(\U_l(\ln)\otimes \onep)=\U_l(\ln)\\
&\label{m2}f(W_l)=char_l(\rho).
\end{align}
\end{lemma}
\begin{proof} We prove claim \ref{m1}: Assume $\omega \otimes \onep\in
  \U_l(\ln)\otimes \onep$. It follows $\omega \in \U_l(\ln)$. It
  follows
\[ f(\omega \otimes \onep)=\omega \oneg =\omega \in \U_l(\ln) \]
and claim \ref{m1} follows. We prove claim \ref{m2}: Assume 
$z=x\otimes w(y-\rho(y)\onep)\in W_l$. It follows $x\in \U_i(\ln)$
$y\in \lp$, and $w(y-\rho(y)\onep)\in \U_j(\lp)$ with $i+j=l$. It
follows
\[ f(z)=xw(y-\rho(y)\oneg)\in char_l(\rho) \]
and claim \ref{m2} follows. The Lemma is proved.
\end{proof}

\begin{theorem} \label{directsum} There is for all $l\geq 1$ an isomorphism
\[ \U_l(\lg)=\U_l(\ln)\oplus char_l(\rho) \]
of vector spaces.
\end{theorem}
\begin{proof} The Theorem follows from Lemma \ref{inclusion},
  Proposition \ref{injective} and the
  fact the natural map $f:\U_l(\ln,\lp)\rightarrow \U_l(\lg)$ is an
  isomorphism of vector spaces.
\end{proof}

In the following we use the notation in \cite{dixmier} Chapter 7.2. Let $P_{++}$ be the
set of dominant weights for $\lg$ and let $\lambda\in \lh^*$ be the
weight with
\[ L(\lambda+\delta)\cong \H^0(X,\O_X(d))^*.\]
Such an element $\lambda$ is uniquely determined since the module
$L(\lambda+\delta)$ is an irreducible finite dimensional $\lg$-module
and there is a one to one correspondence between $P_{++}$ and the set
of irreducible finite dimensional $\lg$-modules. Let $B$ be a basis
for the roots $R$ of $\lg$. It follows $B=L_i-L_{i+1}$ with
$i=1,..,m+n-1$.
Let $v'\in
L(\lambda+\delta)$ be the unique highest weight vector and define two
left ideals $\I'',\I'\subseteq\U(\lg)$ as follows:
\[ \I''=\U(\lg)\nplus+\sum_{h\in \lh}\U(\lg)(h-\lambda(h)), \]
and
\[ \I'=\I''+\sum_{\beta\in B}\U(\nminus)X^{m_\beta}_{-\beta}.\]
Here we let $m_\beta=\lambda(H_\beta)+1$ and $X_{-\beta}$ be a non
zero
element of $\lg^{-\beta}$. It follows by \cite{dixmier} Proposition
7.2.7 the ideal $\I'$ equals the left annihilator ideal $ann(v)$ in $\U(\lg)$ of the
highest weight vector $v=\l^d$. Let $\I'_l=\I'\cap \U_l(\lg)$.

\begin{lemma} \label{ideals}For all $1\leq l < d$ there is an
  equality
 \[ char_l(\rho)=ann_l(v) \]
of filtrations.
\end{lemma}
\begin{proof} Consider the ideal $\I'_l$ for $1\leq l <d$. By
  definition there is an inclusion $char_l(\rho) \subseteq \I'_l$. We prove
  the reverse inclusion. There is an isomorphism
\[ \H^0(X,\O(d))^*\cong L(\lambda+\delta) \]
where $L(\lambda+\delta)$ is the irreducible $\lg$-module with highest
weight $\lambda$. By Lemma \ref{stable} it follows
$\l^d$ has weight $\lambda=d(L_1+\cdots +L_m)$ in the notation of
\cite{fulton}.
Consider the ideal $\I''$:
\[ \I''=\U(\lg)\nplus+\sum_{x\in \lh}\U(\lg)(x-\lambda(x)).\]
It follows $\I''\subseteq char(\rho)$.
Let $\beta_i\in B$ with $\beta_i=L_i-L_{i+1}, 1\leq i\leq m+n-1$.
Let $0\neq E_{ij}\in \lg^{\beta}$ and let
$0\neq H_\beta \in [ \lg^{\beta}, \lg^{-\beta} ]$. One checks 
$ \lambda(H_\beta)+1=1$ if $1\leq i \leq m-1$,
$\lambda(H_\beta)+1=d+1$ if $i= m$ and
$\lambda(H_\beta)+1=1$ if $m+1\leq  i\leq m+n-1$.
Let $K=\sum_{\beta \in B}\U(\lg)X^{m_\beta}_{-\beta}$ and let
$K_l=K\cap \U_l(\lg)$. Let $D$ be the set of integers $i$ with
$i\in \{1,..,m-1,m+1,..,m+n-1\}$. Let $\beta_i=L_i-L_{i+1}$.
It follows 
\[ K_l=\sum_{\beta_i,i\in D}\U_{l-1}(\lg)X^{m_{\beta_i}}_{-\beta_i} \]
and one checks  $K_l\subseteq char_l(\rho)$ and the claim of the
Lemma follows.
\end{proof}

\begin{theorem} \label{main} There is for all $1\leq l < d$ an isomorphism
\[ \U_l(\lg)\cong \U_l(\ln)\oplus ann_l(v) \]
of vector spaces.
\end{theorem}
\begin{proof} This follows from Theorem \ref{directsum} and Lemma \ref{ideals}.
\end{proof}

\begin{corollary} \label{dim} Let $z_1,..,z_D$ be a basis for $\ln$ where
  $D=mn$. It follows the set 
\begin{align} \label{can1}
 \{z_1^{v_1}\cdots z_D^{v_d}(v): 0\leq \sum v_i \leq l \} 
\end{align}
is a basis for $\U_l(\lg)v$ as vector space.
Moreover
\begin{align} \label{can2}
 dim_K(\U_l(\lg)v)=\binom{D+l}{D}.
\end{align}
\end{corollary}
\begin{proof} The natural surjection $\U_l(\lg)\otimes L_v\rightarrow \U_l(\lg)v$
of $K$-vector spaces induce by Theorem \ref{main} an isomorphism
$\U_l(\ln)\otimes L_v \cong \U_l(\lg)v$
of vector spaces. From this isomorphism and the Poincare-Birkhoff-Witt
Theorem  claim \ref{can1} follows. Since
$dim_K(\U_l(\ln))=\binom{D+l}{D}$ claim \ref{can2} follows and the
Corollary is proved.
\end{proof}

It follows we have constructed a basis for $\U_l(\lg)v$ for all $l\geq
1$ and calculated $dim_K(\U_l(\lg)v)$ as a function of $l$.

\section{On  jet bundles and canonical filtrations}

In this section we study the fiber of the jet bundle as $P$-module. We
prove the fiber of the $l$'th jet bundle equals the $l$'th term in the
canonical filtration.

\textbf{Notation:} Let in the following section $W\subseteq V$ be $K$-vector spaces of
dimension $m$ and $m+n$ and let $P\subseteq \SL(V)$ be the parabolic
subgroup fixing $W$. It follows $X`=\SL(V)/P=\gr(m,m+n)$ is the
grassmannian of $m$-planes in $V$. Let $\O_X(d)$ be the line bundle
coming from the Pl\"{u}cker embedding of $X$.

Recall the exact sequence of $P$-modules from Corollary \ref{pmod1}
\begin{align} \label{exactsequence}
0\rightarrow \H^0(X,\m^{l+1}\O_X(d))\rightarrow
\H^0(X,\O_X(d))\rightarrow \Pr^l_X(\O_X(d))(x)\rightarrow 0
\end{align}
when $1\leq l \leq d$.
Dualize sequence \ref{exactsequence} to get the exact sequence
\begin{align} \label{exacts1}
 0\rightarrow \Pr^l_X(\O_X(d))(x)^*\rightarrow
\H^0(X,\O_X(d))^*\rightarrow^\psi \H^0(X,\m^{l+1}\O_X(d))^*\rightarrow
0.
\end{align}
There is by definition an isomorphism
\[ \H^0(X,\O_X(d))^*\cong \Vl \]
where $\lambda=d(L_1+\cdots +L_m)$. The highest weight vector $v\in
\Vl$ is given by $v=l^d$ where $l^d=\sym^d(\wedge^m W)$. We get an
inclusion 
\[ \U_l(\lg)v\subseteq \Vl=\H^0(X,\O_X(d))^*  \]
of $P$-modules.

\begin{theorem} \label{mainp} There is for $1\leq l <d$ an isomorphism
\[\Pr^l_X(\O_X(d))(x)^*\cong \U_l(\lg)v \]
of $P$-modules.
\end{theorem}
\begin{proof}  Consider the exact sequence \ref{exacts1} of $P$-modules
\[ 0\rightarrow \Pr^l_X(\O_X(d))(x)^*\rightarrow
\H^0(X,\O_X(d))^*\rightarrow^\psi \H^0(X,\m^{l+1}\O_X(d))^*\rightarrow
0.\]
There is an inclusion of $P$-modules 
\[ \U_l(\lg)v\subseteq \H^0(X,\O_X(d))^* \]
where $v$ is the highest weight vector.
Consider an element $x_1\cdots x_iv\in \U_l(\lg)v$ with $i\leq l$ and
  $x_i\in  \lg$.  It follows
\[ \psi(x_1\cdots x_iv)(s)=x_1\cdots x_is(e) \]
for $s\in \H^0(X,\m^{l+1}\O_X(d))$. The section $s$ has a zero of
order $\geq l+1$ at $e$. Since $x_1\cdots x_i$ acts as a differential
operator of order $i$ it follows
\[ x_1\cdots x_is\in \H^0(X,\m^{l+1-i}\O_X(d)) \]
hence $x_1\cdots x_is$ has a zero of order $l+1-i$ at $e$. It follows 
$\psi(x_1\cdots x_iv)=0$ since $i\leq l$. We get $\psi(\U_l(\lg)v)=0$ 
and $\U_l(\lg)v\subseteq ker(\psi)=\Pr^l_X(\O_X(d))(x)^*$ since the
sequence above is an exact sequence of $P$-modules.
We get an inclusion
\[ \U_l(\lg)v\subseteq \Pr^l_X(\O(d))(x)^* \]
of $P$-modules when $1\leq l <d$.
By Corollary \ref{dim} it follows 
\[ dim_K(\U_l(\lg)v)=\binom{mn+l}{mn}=dim_K(\Pr^l_X(\O_x(d))(x)^*\]
hence the Theorem is proved.
\end{proof}

Note: In \cite{kallstrom}, Section 5 the authors claim they prove a similar
result using different techniques. The aim of the introduction of the
techniques in this paper is to use them to get more precise
information on $\U_l(\lg)v$ as $P$-module. There is work in progress
(see \cite{maa31}) on giving a more detailed description of the
$P$-module $\U_l(\lg)v$ and to apply this to the study of resolutions
of discriminants of linear systems on grassmannians.

Let $\underline{d}=(d_1,..,d_k)$ with $d_i\geq 1$. 
Let $\E(\underline{d})=\O_X(d_1)\oplus \cdots \oplus
\O_X(d_k)$ with $X=\gr(m,m+n)$. Let
$\H^0(X,\O_X(d_i))^*=V_{\lambda_i}=V_i$ with highest weight vector
$v_i$. Let $W=\{v_1,..,v_k\}\subseteq \H^0(X,\E(\underline{d}))^*$ be
the subspace spanned by the vectors $v_i$.
\begin{corollary} There is for $1\leq l \leq min\{d_i\}$ an
  isomorphism
\[ \Pr^l_X(\E(\underline{d}))(x)^*\cong
\U_l(\lg)W=\oplus_{i=1}^k\U_l(\lg)v_i \]
of $P$-modules.
\end{corollary}
\begin{proof}
We get from Theorem \ref{mainp} the following:
\[ \Pr^l_X(\E(\underline{d}))(x)^*\cong
\oplus_{i=1}^k\Pr^l_X(\O_X(d_i))(x)^*\cong
\oplus_{i=1}^k\U_l(\lg)v_i\cong \U_l(\lg)W\]
and the Corollary follows.
\end{proof}

\section{Discriminants of linear systems on the grassmannian}

In this section we use the results obatined in the previous sections
to prove the $l$'th discriminant of any linear system on any
grassmannian is irreducible.

\textbf{Notation:} Let in the following section $X=\gr(m,m+n)$ be the grassmannian of $m$-planes in an
$m+n$-dimensional vector space and let $\O_X(1)$ be the line bundle
coming from the Plucker embedding. Let $\O_X(d)=\O_X(1)^{\otimes
  d}$. Let $S=\Spec(K)$ be the spectrum of the base field $K$.

Let $1\leq l \leq d$ and consider the exact sequence \ref{proj2} of locally free sheaves
\begin{align}
&\label{grassex} 0\rightarrow \Q_{l,d}\rightarrow \H^0(X,\O_X(d))\otimes \O_X \rightarrow \Pr^l_X(\O_X(d)) \rightarrow 0 
\end{align}
where
\[ \Q_{l,d}=p_*(\I_{\Delta}^{l+1}\otimes q^*\O_X(d)).\]
It follows from \cite{maa10}, Theorem 2.5 we get a commutative diagram
\[
\diagram  \p(\Q_{l,d}^*)\rto^i \dto^\pi & \p(W^*)\times_S X \rto^p \dto^q
& X \dto \\
          D^l(\O_X(d)) \rto^j & \p(W^*) \rto & S
\enddiagram
\]
where $W=\H^0(X,\O_X(d))$. Here $i,j$ are closed immersions of schemes
and $\pi$ is the restriction of $q$. It followd $D^l(\O_X(d))$ is the
scheme theoretic direct image of $\p(\Q_{l,d}^*)$. 

\begin{theorem} \label{irreducible} The discriminant $D^l(\O_X(d))$ is irreducible when
  $1\leq l \leq d$.
\end{theorem}
\begin{proof}
Since $\Q_{l,d}^*$ is locally free it follows from \cite{maa10},
Corollary 2.6 $D^l(\O_X(d))$ is irreducible.
\end{proof}

In a series of papers (see \cite{maa1}, \cite{maa10}, \cite{flag},
\cite{maa2}, \cite{maa3} and \cite{maa4}) the structure of the jet
bundle $\Pr^l_X(\O_X(d))$ as left and right $\O_X$-module and left and
right $P$-module has been studied using various techniques: Explicit
techniques, group theoretic techniques and Lie theoretic
techniques. This study is part of a project where the aim is to study
discriminants of linear systems on flag varieties (see \cite{maa10}
and \cite{flag}). 
The aim of the study of
$\Pr^l_{X}(\O_X(d))$ is to use its properties to study the map $\pi$ and $D^l(\O_X(d))$.
We want information on the syzygies of $\D^l(\O_X(d))$, the singularity type
of $D^l(\O_X(d))$, its degree and its dimension.

\begin{example} Syzygies of discriminants of linear systems on grassmannians.\end{example}

By \cite{maa09} the following holds: Let $V=\H^0(X,\O_X(d))$ and 
consider the exact sequence \ref{grassex} of locally free
sheaves. Dualize it to get the short exact sequence

\[ 0 \rightarrow \Pr^l_X(\O_X(d))^*\rightarrow V^*\otimes_{\O_S} \O_X
\rightarrow \Q_{l,d}^*\rightarrow 0.\]
Take relative projective space bundle to get the closed immersion of
schemes
\[ \p(\Q_{l,d}^*)\rightarrow^i \p(V^*\otimes_{\O_S} \O_X)\cong \p(V^*)\times_S X.\]
Consider the sequence of locally free sheaves on $Y=\p(V^*)\times_S X$:
\[ \O(-1)_Y\rightarrow V\otimes \O_Y\rightarrow^{T^l_Y}
\Pr^l_X(\O_X(d))_Y .\]

Let the composed morphism be $\phi_{l,d}$.
By definition (see \cite{maa09} and \cite{maa10}) the zero scheme
$Z(\phi_{l,d})$ equals the first incidence scheme $I_1(T^l)$ of the
$l$'th Taylor morphism. It also follows from \cite{maa10} there is an equality of subschemes
\[ \p(\Q_{l,d}^*)=I_1(T^l)\]
of $\p(V^*)\times_S X$
The direct image scheme $p(\p(Q_{l,d}^*))$ equals the $l$'th
discriminant $D_l(\O_X(d))$ of the line bundle $\O_X(d)$. Since
$I_1(\O_X(d))$ is defined 
as the zero scheme of the map $\phi_{l,d}$ of locally free sheaves
there is by \cite{maa09} a \emph{Koszul complex} of locally free sheaves
\begin{align}
&\label{koszul} 0\rightarrow \O_X(-r)_Y\otimes \wedge^r
\Pr^l_X(\O_X(d))_Y^*\rightarrow \cdots \rightarrow \O_X(-2)_Y\otimes
\wedge^2 \Pr^l_X(\O_X(d))_Y^*\rightarrow 
\end{align}
\[ \O_X(-1)_Y\otimes \Pr^l_X(\O_X(d))_Y^*\rightarrow \O_Y \rightarrow
\O_{I_1(T^l)}\rightarrow 0 \]
of $\O_Y$-modules.
Since $X$ is smooth it is Cohen-Macaulay and it follows from \cite{maa09} that \ref{koszul} is
a resolution of $\O_{I_1(T^l)}$. When we push the complex
\ref{koszul} down to $\p(V^*)$ we get
a double complex with terms given as follows:
\begin{align}
&\label{double}  C^{i,j}(T^l)=\O_X(-i)\otimes_S \H^j(X,\wedge^i
\Pr^l_X(\O_X(d))^*).
\end{align}
It is hoped the couble complex \ref{double} will give information on
the syzygies of the discriminant $D^l(\O_X(d))$ for all  $1\leq l <d$.
One needs to calculate the cohomology groups $\H^j(X,\wedge^i
\Pr^l_X(\O_X(d))^*)$ for all $i,j$. The aim of the study done in this
paper is to use  information on the canonical filtration studied in the previous
sections to get such information and to decide if $C^{i,j}(T^l)$ may
be used. Note: By \cite{maa10} the discriminant $D_1(\O(d))$ on $\p^1$
is a determinantal scheme. For determinantal schemes much is known
about their syzygies due to the work of Lascoux (see
\cite{lascoux}). It follows we get two approaches to the study of
syzygies of discriminants of linear systems on flag schemes.
 
There is by the previous sections an isomorphism
\[ \Pr^l_X(\O_X(d))^* \cong \U^l(\lg)v \]
of $P$-modules for all $1\leq l <d$. We get for all $1\leq i \leq rk(\Pr^l_X(\O_X(d)))$ 
an isomorphism
\[ \wedge^i \Pr^l_X(\O_X(d))(x)^* \cong \wedge^i \U^l(\lg)v \]
of $P$-modules.
There is a filtration
\[ 0 \neq U_1\subsetneq U_2 \subsetneq \cdots \subsetneq U_k=\wedge^i
\U^l(\lg)v \]
of $P$-modules with $U_n/U_{n-1}$ an irreducible $P$-module for all $n=2,..,k$. Let $\E_n$ be the
locally free $\O_X$-module corresponding to $U_n/U_{n-1}$. It follows
$\E_n$ is an irreducible $\O_X$-module with a $G$-linearization where
$G/P=\gr(m,m+n)$. It follows the cohomology groups $\H^j(X,\E_n)$
are determined by Bott's Theorem. This approach may be used to
calculate the cohomology group $\H^j(X,\wedge^i \Pr^l_X(\O_X(d))^*)$
for all $i,j$. We aim to get a complete description of the double complex
$C^{i,j}(T^l)$ using higher direct images of sheaves, annihilator
ideals of highest weight vectors, Taylor maps and Bott's Theorem.

\textbf{Acknowledgements}. The author  thanks Michel Brion, Alexei Roudakov and Claire Voisin for 
interesting discussions on subjects related to the problems discussed in this paper.

\end{document}